\documentclass[12pt]{article}
\usepackage{mathrsfs}
\usepackage{amsfonts}
\usepackage[all]{xy}
\usepackage{amsmath,amssymb}
\usepackage{amscd}

\baselineskip 5.5pt \lineskip 5.5pt \numberwithin{equation}{section}

\def\qed{\hfill$\Box$\par}

\def\qed{\ \ \ifhmode\unskip\nobreak\fi\ifmmode\ifinner
         \else\hskip5pt\fi\fi
 \hbox{\hskip5pt\vrule width4pt height6pt depth1.5pt\hskip 1 pt}}

\begin{document}
\title{\large{\textbf{Tensor product weight representations of\\ the
Neveu-Schwarz algebra}} }
\author{Xiufu Zhang\\
{\scriptsize  School of Mathematics and Statistics, Jiangsu Normal
University, Xuzhou 221116, China}}
\date{}
\maketitle \footnotetext{\footnotesize Supported partially by the
National Natural Science Foundation of China (No. 11271165) and the
Youth Foundation of National Natural Science Foundation of China
(No.11101350).} \footnotetext{\footnotesize* Email:
xfzhang@jsnu.edu.cn}

\numberwithin{equation}{section}

\begin{abstract}
In this paper, the tensor product of highest weight modules with intermediate series modules over the Neveu-Schwarz algebra is studied. The weight spaces of  such tensor products are all infinitely dimensional if the highest weight module is nontrivial. We find that all such tensor products are indecomposable. We give the necessary and sufficient conditions for these tensor product modules to be irreducible by using ``shifting technique" established for the Virasoro case in [13]. The necessary and sufficient conditions for any two such tensor products to be isomorphic are also determined.
 \vspace{2mm}\\{\bf 2000 Mathematics Subject
Classification:} 17B10, 17B65, 17B68 \vspace{2mm}
\\ {\bf Keywords:}  Neveu-Schwarz algebra, Tensor product module, Indecomposable module,
 Irreducible module, Isomorphism.
\end{abstract}

\vskip 3mm\noindent{\section{Introduction}}

\vskip 3mm We denote by $\mathbb{Z}, \mathbb{Z}_+, \mathbb{N}$ and
$\mathbb{C}$ the set of all integers, nonnegative integers, positive
integers and complex numbers, respectively. In this paper, we set
$Z=\frac{1}{2}+\mathbb{Z}.$

A Lie superalgebra is a generalization of a Lie algebra to include a $\mathbb{Z}_2$-grading. Lie superalgebras are important in theoretical physics where they are used to describe the mathematics of supersymmetry (see [1], [2] and the references there in).
A super Virasoro algebra is an extension of the Virasoro algebra to a Lie superalgebra. There are two extensions with particular importance in superstring theory: the Ramond algebra and  the Neveu-Schwarz algebra. They describe the symmetries of a superstring in two different sectors, called the Ramond sector and the Neveu-Schwarz sector. Concretely, the
super-Virasoro algebra  is the Lie superalgebra
$SV(\theta)=SV(\theta)_{\bar{0}}\oplus SV(\theta)_{\bar{1,}}$ where
$\theta=\frac{1}{2}$ or $0$ according as $SV(\theta)$ is the
Neveu-Schwarz algebra $\mathfrak{NS}$ or the Ramond algebra $\mathfrak{R}$, $SV(\theta)_{\bar{0}}$
has a basis $\{L_n, C\mid n\in \mathbb{Z}\}$ and
$SV(\theta)_{\bar{1}}$ has a basis $\{G_r\mid r\in
\theta+\mathbb{Z}\},$ with the commutation relations
\begin{align*}
&[L_m,L_n]=(n-m)L_{n+m}+\delta_{m+n,0}\frac{m^{3}-m}{12}C,\\
&[L_m,G_r]=(r-\frac{m}{2})G_{m+r},\\
&[G_r,G_s]=2L_{r+s}-\frac{1}{3}\delta_{r+s,0}(r^2-\frac{1}{4})C,\\
&[SV_{\bar{0}},C]=0=[SV_{\bar{1}},C],
\end{align*}
for $m, n\in\mathbb{Z}, r, s \in \theta+\mathbb{Z}.$ It is obvious
that the Neveu-Schwarz algebra is $\frac{1}{2}\mathbb{Z}$-graded and
Ramond algebra is $\mathbb{Z}$-graded:
$\mathrm{deg}L_i=i,\mathrm{deg}G_r=r$ and $\mathrm{deg}C=0,$ which
induces the grading of the universal algebra $U(SV(\theta))$ of $SV(\theta).$

The representation theory over super-virasoro algebras has been attracting a lot of attentions, see for example [3]-[8] and the references there
in. There are two classic families of simple Harish-Chandra modules over the super-Virasoro algebras: highest (lowest) weight modules and the  intermediate series modules. It is proved in [5] that these two families modules exhaust all simple Harish-Chandra modules.  In [8] it is even shown that the simple weight modules with a finite-dimensional weight space for $\mathfrak{NS}$ are all Harish-Chandra modules.

Several classes of modules over the Virasoro algebra are constructed from tensor products of simple Virasoro modules. For details, we refer the readers to [9], [10], [11]. According to [12], all weight spaces of such simple modules are infinite-dimensional. The irreducibility problem for the tensor products over the Virasoro algebra is an interesting problem. In [13], the irreducibility for the tensor product of highest weight modules with intermediate series modules over the Virasoro algebra was determined by using Feigin-Fuchs' Theorem (see [14]) on singular vectors of Verma modules and the so called ``shifting technique".

In this paper, we will study the tensor product modules over the Neveu-Schwarz algebra $\mathfrak{NS}.$
Let us denote by $\mathfrak{h}$ the Cartan subalgebra of $\mathfrak{NS}$ with
the basis $L_0$ and $C$, by $\mathfrak{n}^{-}$ the Lie subalgebra
with the basis $\{E_{-\frac{i}{2}}|i\in\mathbb{N}\}$ and by
$\mathfrak{n}^{+}$ the Lie subalgebra with the basis
$\{E_{\frac{i}{2}}|i\in\mathbb{N}\}$,
where, and in the following, $$E_s=\left\{\begin{array}
{l@{\quad\quad}l} L_s, & \mathrm{if}\ \ s\in \mathbb{Z}, \\
G_s, &\mathrm{if}\ \ s\in Z.
\end{array}\right.$$
We also denote by $\mathfrak{b}^+$ the Lie subalgebra with the basis
$\{L_i,G_{\frac{1}{2}+i},C|i\in\mathbb{Z_+}\}.$

In [5], a class of representations $SA_{a,b}$ over $\mathfrak{NS}$
for $a,b\in\mathbb{C}$ are defined. Concretely,  $SA_{a,b}=(SA_{a,b})_{\bar{0}}\oplus (SA_{a,b})_{\bar{1}},$ where  $(SA_{a,b})_{\bar{0}}$ has a basis $\{x_i|i\in\mathbb{Z}\},$ and $(SA_{a,b})_{\bar{1}}$ has a basis $\{y_{\frac{1}{2}+i}|i\in\mathbb{Z}\},$ the action of $\mathfrak{NS}$ on  $SA_{a,b}$  is given
by
\begin{align*}
&L_ix_j=(a+j+ib)x_{i+j},\\
&L_iy_{\frac{1}{2}+j}=(a+\frac{1}{2}+j+i(b-\frac{1}{2}))y_{\frac{1}{2}+i+j},\\
&G_{\frac{1}{2}+i}x_{j}=y_{\frac{1}{2}+i+j},\\
&G_{\frac{1}{2}+i}y_{\frac{1}{2}+j}=(a+\frac{1}{2}+j+2(\frac{1}{2}+i)(b-\frac{1}{2}))x_{1+i+j},\\
&Cx_i=Cy_{\frac{1}{2}+j}=0.
\end{align*}
The following theorem is taken from [5] and it can be checked directly:
\vskip 3mm\noindent{\bf{Theorem 0.1.}}
(1) $SA_{a,b}\simeq SA_{c,d}$ if and only if one of the following holds:
$$a-c\in\mathbb{Z},b=d;$$
$$a\notin
\mathbb{Z},a-c\in\frac{1}{2}+\mathbb{Z},b=1,d=\frac{1}{2};$$
$$a\notin\frac{1}{2}+\mathbb{Z},a-c\in\frac{1}{2}+\mathbb{Z},b=\frac{1}{2},d=1.$$

(2) If $a\in \mathbb{Z},$ then $SA_{a,1}\cong SA_{0,1}$ and $SA_{0,1}$ has a unique simple submodule $SA_{0,1}^{'}$ spanned by
$\{x_j|j\in\mathbb{Z},j\neq
0\}\cup\{y_{\frac{1}{2}+j}|j\in\mathbb{Z}\}.$
If $a\in\frac{1}{2}+\mathbb{Z},$ then $SA_{a,\frac{1}{2}}\cong SA_{\frac{1}{2},\frac{1}{2}}$ and
$SA_{\frac{1}{2},\frac{1}{2}}^{'}=SA_{\frac{1}{2},\frac{1}{2}}/\mathbb{C}y_{-\frac{1}{2}-n}$
is irreducible.

Let $SA_{a,b}^{'}=SA_{a,b}$ if $SA_{a,b}^{'}$ is irreducible and let $SA_{0,1}^{'}, SA_{\frac{1}{2},\frac{1}{2}}^{'}$ be defined as above. Then we get a class of irreducible modules
$\{SV_{a,b}^{'}|a,b\in\mathbb{C}\},$ which is called the
indeterminate series over $\mathfrak{NS}$.

As a principal observation, we have the following theorem:

\vskip 3mm\noindent{\bf{Theorem 0.2.}} $SA_{0,1}^{'}\cong SA_{-\frac{1}{2},\frac{1}{2}}^{'}\cong SA_{\frac{1}{2},\frac{1}{2}}^{'}.$

In fact, it is easy to check that\begin{align*}
f:SA_{0,1}^{'} \rightarrow SA_{-\frac{1}{2},\frac{1}{2}}:
 x_i  \mapsto \frac{1}{i}y_{\frac{1}{2}+i},
 y_{\frac{1}{2}+j}  \mapsto x_{j+1},
\end{align*}
is an isomorphism.

By Theorems 0.1 and 0.2, for $SA_{a,b}^{'}$  we can assume that $0\leq\mathfrak{Re}a<1, b\neq1$ in this paper.

For any $h,c\in\mathbb{C},$ we can view a 1-dimensional vector space
$\mathbb{C}u$ as a $b^+$-module such that $\mathfrak{n}^+u=0, L_0u=hu$ and
$Cu=cu.$ Then we obtain the Verma module over $\mathfrak{NS}$:
$$M(c,h):=U(\mathfrak{NS})\bigotimes_{U(b^{+})}\mathbb{C}u.$$

It is well known that the Verma module $M(c,h)$ has a unique maximal
submodule $J(c,h)$ and the corresponding simple quotient module is
denoted by $V(c,h).$ A nonzero weight vector $u^{'}\in M(c,h)$ is
called a singular vector if $\mathfrak{n}^{+}u^{'}=0.$ It is clear that
$J(c,h)$ is generated by all singular vectors in $M(c,h)\setminus
\mathbb{C}u$ and that $M(c,h)=V(c,h)$ if and only if there does not
exist any other singular vectors besides $\mathbb{C}u$.  The
singular vectors in $M(c,h)$ are described explicitly in [3]. In
particular,  $J(c,h)$ can be generated by at most two singular
vectors. If $J(c,h)$ is generated by two singular vectors, we can
find homogeneous $Q_1, Q_2\in U(\mathfrak{n}^-)$ such that
$$J(c,h)=U(\mathfrak{n}^-)Q_1u+U(\mathfrak{n}^-)Q_2u.$$
If $J(c,h)$ is generated by one singular vectors, we can find the
unique $Q_1\in U(\mathfrak{n}^-)$ up to a scalar multiple such that
$J(c,h)=U(\mathfrak{n}^-)Q_1u,$ in which case we set $Q_2=Q_1$ for convenience.
When $M(c,h)=V(c,h)$ we set $Q_1=Q_2=0.$

In this paper, we mainly study the tensor product of simple highest weight modules with intermediate series modules over the Neveu-Schwarz algebra $\mathfrak{NS}$.
In section 2, we prove that the tensor product module $V\otimes SA_{a,b}^{'}$ of any nontrivial highest weight module $V$ with $SA_{a,b}^{'}$ are indecomposable. Moreover, if $V=M(c,h)$ is a Verma module, then $V\otimes SA_{a,b}^{'}$  is reducible.
In section 3, we determine the necessary and sufficient conditions for $V(c,h)\otimes SA_{a,b}^{'}$ to be simple.
In section 4, we obtain that any two tensor products of the form $V(c,h)\otimes SA_{a,b}^{'}$ are isomorphic if and only if the corresponding highest weight modules and intermediate series modules are isomorphic respectively.

\vskip 3mm\noindent{\section{Indecomposable objects}}

 Let $V$ be a highest weight module (not necessarily simple) over $NS$ with highest
weight vector $u$ of  weight $(c,h).$ Since $V$ and $SA_{a,b}^{'}$ are
$L_0$-diagonalizable, so is $V\otimes SA_{a,b}^{'}:$
$$V\otimes SA_{a,b}^{'}=\bigoplus_{m\in\mathbb{Z}\cup Z}(V\otimes SA_{a,b}^{'})_{m+h+a},$$
where
$$(V\otimes SA_{a,b}^{'})_{m+h+a}=\bigoplus_{n\in\mathbb{Z}\cup Z}V_{h+n}\otimes
\mathbb{C}v_{m-n},$$ and
$$v_{m-n}=\left\{\begin{array}
{l@{\quad\quad}l} x_{m-n}, & \mathrm{if}\ m-n\in \mathbb{Z}, \\
y_{m-n},  &\mathrm{if} \ m-n\in\frac{1}{2}+\mathbb{Z}.
\end{array}\right.$$

\vskip 3mm\noindent$\mathbf{Remark}$:  If $V$ is nontrivial, then
$(V\otimes SA_{a,b}^{'})_{m+h+a}$ is infinite dimensional for all
$m\in \mathbb{Z}\cup Z$.

\vskip 3mm\noindent{\bf{Lemma 2.1.}} The module $V\otimes SA_{a,b}^{'}$ is
generated by  $\{u\otimes v_m|m\in\mathbb{Z}\cup Z\}.$

\vskip 3mm\noindent{\bf{Proof.}} Note that $V\otimes SA_{a,b}^{'}$
is spanned by $\{wu\otimes v_m|m\in \mathbb{Z}\cup Z, w\in
U(\mathfrak{n}^-)\},$ so the lemma holds. \hfill$\Box$

\vskip 3mm The following theorem shows that $V\otimes SA_{a,b}^{'}$
is an indecomposable $\mathfrak{NS}$-module.

\vskip 3mm\noindent{\bf{Theorem 2.2.}} Let $a,b,c,h\in \mathbb{C}$, $V$ be any highest weight
module with highest weight vector $u$ of weight $(c,h)$ over the Neveu-Schwarz algebra. Then
$$\mathrm{End}_{\mathfrak{NS}}(V\otimes SA_{a,b}^{'})\cong\mathbb{C}.$$

\vskip 3mm\noindent{\bf{Proof.}} For any
$\psi\in\mathrm{End}_{\mathfrak{NS}}(V\otimes SA_{a,b}^{'}),$ $\psi(u\otimes
v_m)$ and $u\otimes v_m$ have the same weight, so
$$\psi(u\otimes v_m)=u_0\otimes v_m+u_{\frac{1}{2}}\otimes
v_{m+\frac{1}{2}}+u_{1}\otimes v_{m+1}+\cdots+u_{\frac{n}{2}}\otimes
v_{m+\frac{n}{2}},$$ where $u_{\frac{i}{2}}\in
V_{h-\frac{i}{2}},i=0,1,\cdots,n, u_{\frac{n}{2}}\neq0.$

\vskip 3mm\noindent{\bf{Claim 1.}} $n=0.$

In fact, when $m\in Z,$ i.e.,
$v_m=y_m,$ we see that $G_{\frac{1}{2}+l}y_m\neq0$ for sufficiently
large $l.$ Otherwise, $a+m+(2l+1)(b-\frac{1}{2})=0$ for infinitely
many $l$, which means that $b=\frac{1}{2}$ and $a=-m\in Z.$
By Theorem 0.1 we can assume that $a=\frac{1}{2},$ so
$y_{-\frac{1}{2}}\in SA_{a,\frac{1}{2}}^{'}$
 which contradicts the definition of $SA_{\frac{1}{2},\frac{1}{2}}^{'}$.
For sufficiently large $l,$ we define
$$w=L_{2l+2}-\frac{a+m+(2l+2)(b-\frac{1}{2})}{a+m+(2l+1)(b-\frac{1}{2})}G_{l+\frac{3}{2}}G_{l+\frac{1}{2}}.$$
Since $wu=0,$ we have $$w\psi(u\otimes y_m)=0,$$ i.e.,
$$u_{\frac{n}{2}}\otimes wv_{m+\frac{n}{2}}=0$$ since
$wu_{\frac{i}{2}}=0, i=0,1,\cdots,n$ for sufficiently large $l.$
If $n$ is odd, then
$$
a+m+\frac{n}{2}+(2l+2)b-\frac{a+m+(2l+2)(b-\frac{1}{2})}{a+m+(2l+1)(b-\frac{1}{2})}(a+m+\frac{n}{2}+l+\frac{1}{2}+(2l+3)(b-\frac{1}{2}))=0
$$for sufficiently large $l,$
so the coefficients of $l$ are all zero, and then we can deduce that
$b=\frac{1}{2},a=-m\in Z,$ which is also absurd by the definition of
$SA_{a,b}^{'}$. So $n$ is a even number, and
$$
a+m+\frac{n}{2}+(2l+2)(b-\frac{1}{2})-\frac{a+m+(2l+2)(b-\frac{1}{2})}{a+m+(2l+1)(b-\frac{1}{2})}(a+m+\frac{n}{2}+(2l+1)(b-\frac{1}{2}))=0
$$for sufficiently large $l.$ Moreover, we can deduce that
$b=\frac{1}{2}$ or $n=0.$ If $n\neq0,$ then $b=\frac{1}{2}.$ For
sufficiently large $l$, $L_l^{2}y_m\neq0.$ Otherwise, we can deduce
that $a=-m,$ which is again absurd. For sufficiently large $l$, set
$$w^{'}=L_{2l}-\frac{a+m}{(a+m)(a+m+l)}L_l^2.$$ Then
$w^{'}\psi(u\otimes y_m)=0,$  furthermore,
$w^{'}y_{m+\frac{n}{2}}=0.$ So we can deduce that
$$\frac{n}{2}(a+m)(a+m+\frac{n}{2})=0,$$ and we see that
$a=-m-\frac{n}{2}\in Z,$ which is impossible by the definition of
$SA_{a,b}^{'}$. Thus
$b\neq\frac{1}{2}$ and $n=0.$

When $m\in \mathbb{Z},$ i.e., $v_m=x_m,$ we see that
$L_lG_{\frac{1}{2}}x_m\neq0$ for sufficiently large $l.$ Otherwise,
$a+\frac{1}{2}+m+l(b-\frac{1}{2})=0$ for infinitely many $l$, which
means that $b=\frac{1}{2}$ and $a=-\frac{1}{2}-m\in Z,$ a contradiction. For sufficiently
large $l,$ we define
$$w=G_{l+\frac{1}{2}}-\frac{1}{a+m+\frac{1}{2}+l(b-\frac{1}{2})}L_{l}G_{\frac{1}{2}}.$$
Since $wu=0,$ we have $$w\psi(u\otimes x_m)=0,$$ i.e.,
$$u_{\frac{n}{2}}\otimes wv_{m+\frac{n}{2}}=0$$ since
$wu_{\frac{i}{2}}=0, i=0,1,\cdots,n$ for sufficiently large $l.$
If $n$ is odd, then
$$a+m+\frac{n}{2}+(2l+1)(b-\frac{1}{2})-\frac{(a+m+\frac{n}{2}+(b-\frac{1}{2}))(a+m+\frac{n}{2}+\frac{1}{2}+lb)}{a+m+\frac{1}{2}+l(b-\frac{1}{2})}=0$$
for sufficiently large $l.$ We can deduce that $b=\frac{1}{2},
a=-m-\frac{n}{2}\in Z,$ a contradiction again. Hence $n$ is even and
$$1-\frac{a+\frac{1}{2}+m+\frac{n}{2}+l(b-\frac{1}{2})}{a+m+\frac{1}{2}+l(b-\frac{1}{2})}=0.$$
So we have $n=0.$ This completes the proof of Claim 1.

Now we know that for any $m\in\mathbb{Z}\cup Z,$ there exists a
$c_m\in \mathbb{C}$ such that
$$\psi(u\otimes v_m)=c_m(u\otimes v_m).$$
For $m\in \mathbb{Z},$ if $L_{l}y_{\frac{1}{2}+m}=0$ for any
sufficiently large $l\in\mathbb{Z},$ then
$a+\frac{1}{2}+m+l(b-\frac{1}{2})=0,$ which means $b=\frac{1}{2},a=-\frac{1}{2}-m\in
Z,$   which derives a contradiction. So we can choose a sufficiently
large $l\in \mathbb{N}$ such that
$$
L_{l}(u\otimes y_{m+\frac{1}{2}})\neq0.$$ By
\begin{eqnarray*}
G_{l+\frac{1}{2}}\psi(u\otimes x_m)&=&c_m(u\otimes y_{\frac{1}{2}+l+m})\\
&=&c_{\frac{1}{2}+m+l}(u\otimes y_{\frac{1}{2}+l+m}),
\end{eqnarray*}
and
\begin{eqnarray*}
L_{l}\psi(u\otimes y_{\frac{1}{2}+m})&=&c_{m+\frac{1}{2}}(a+m+\frac{1}{2}+l(b-\frac{1}{2}))(u\otimes y_{\frac{1}{2}+l+m})\\
&=&c_{\frac{1}{2}+m+l}(a+m+\frac{1}{2}+l(b-\frac{1}{2}))(u\otimes
y_{\frac{1}{2}+l+m}),
\end{eqnarray*}
we obtain that $c_m=c_{m+\frac{1}{2}}=c_{\frac{1}{2}+m+l},\ \
m\in\mathbb{Z}.$ Similarly, for $m\in Z,$ we can also obtain that
$c_m=c_{m+\frac{1}{2}}.$ Thus
$$c_m=c_{m+\frac{1}{2}}=c_{\frac{1}{2}+m+l},\ \ m\in\mathbb{Z}\cup Z.$$
Therefore, $\psi(u\otimes v_m)=c(u\otimes v_m),$ where $c$ is a
constant. By Lemma 2.1, $\psi$ is a scalar. Thus
$$\mathrm{End}_{\mathfrak{NS}}(V\otimes SA_{a,b}^{'})\cong \mathbb{C}.$$
This completes the proof of Theorem 2.2.\hfill$\Box$

\vskip 3mm The following theorem shows that $M(c,h)\otimes
SA_{a,b}^{'}$ is always reducible, even if $M(c,h)$ is irreducible.

\vskip 3mm\noindent{\bf{Theorem 2.3.}} $M(c,h)\otimes SA_{a,b}^{'}$
is reducible for all $a,b,c,h\in\mathbb{C}.$

\vskip 3mm\noindent{\bf{Proof.}} It is sufficient to prove that
every $u\otimes v_m$ generates a proper submodule of $M(c,h)\otimes
SA_{a,b}^{'}$, where $u$ is the highest weight vector of $M(c,h)$.
Assume that $M(c,h)\otimes SA_{a,b}^{'}$ is cyclic on $u\otimes v_m,$
i.e., $$M(c,h)\otimes SA_{a,b}^{'}=U(\mathfrak{NS})(u\otimes
v_m)=U(\mathfrak{n}^{-})U\mathfrak{(n}^{+})(u\otimes v_m).$$ Then there must exists a
$w\in U(\mathfrak{n}^{-})U(\mathfrak{n}^{+})$ such that $$u\otimes
v_{m-\frac{1}{2}}=w(u\otimes v_m).$$ Let $$w=\sum_{\begin{subarray}\
k_1,k_2,\cdots,k_n \in \mathbb{Z}_{+}, \\
l_i=0\ \mathrm{or}\  1,1\leq i\leq r
\end{subarray}}L_{-n}^{k_n}\cdots\L_{-1}^{k_1}G_{-\frac{1}{2}-r}^{l_r}\cdots
G_{-\frac{1}{2}-1}^{l_1}G_{-\frac{1}{2}}^{l_0}x_{k_1,\cdots,k_n}^{l_0,\cdots,l_r},
$$
where $x_{k_1,\cdots,k_n}^{l_0,\cdots,l_r}\in U(\mathfrak{n}^{+})$ is a
homogeneous element. Since $$x_{k_1,\cdots,k_n}^{l_0,\cdots,l_r}u=0,\
\forall\ x_{k_1,\cdots,k_n}^{l_0,\cdots,l_r}\in
U(\mathfrak{n}^+)\setminus\mathbb{C},$$ we can assume there exists some
$x_{k_1,\cdots,k_n}^{l_0,\cdots,l_r}v_m\neq0.$ Let
$$L_{-n}^{k_n}\cdots\L_{-1}^{k_1}G_{-\frac{1}{2}-r}^{l_r}\cdots
G_{-\frac{1}{2}-1}^{l_1}G_{-\frac{1}{2}}^{l_0}$$ be a term in the
expression of $w$ such that
$\sum_{i}ik_i+\sum_{j}(\frac{1}{2}+j)l_j$ is maximal. By comparing two
sides of $u\otimes v_{m-\frac{1}{2}}=w(u\otimes v_m),$ we have
$$L_{-n}^{k_n}\cdots\L_{-1}^{k_1}G_{-\frac{1}{2}-r}^{l_r}\cdots
G_{-\frac{1}{2}-1}^{l_1}G_{-\frac{1}{2}}^{l_0}u\otimes
x_{k_1,\cdots,k_n}^{l_0,\cdots,l_r}v_m=0.$$ Since $M(c,h)$ is a free
$U(\mathfrak{n}^{-})$-module, it follows that
$x_{k_1,\cdots,k_n}^{l_0,\cdots,l_r}v_m=0,$
which is a contradiction. This completes the proof of Theorem 2.3.\hfill$\Box$

\vskip 3mm\noindent{\section{Irreducible objects}}

For $a, b, c, h\in \mathbb{C}$ with $0\leq \mathfrak{Re}a<1$ and
$b\neq 1,$
by Theorem 2.3, $M(c,h)\otimes
SA_{a,b}^{'}$ is always reducible, even if $M(c,h)$ is irreducible.
  We will give the necessary and sufficient conditions for the simplicity of $V(c,h)\otimes SA_{a,b}^{'}$ in this section.

Recall that $U(\mathfrak{NS})$ has an natural $\mathbb{Z}\cup Z$-gradation and an induced $\mathbb{Z}_2$-gradation. For homogeneous $p\in U(\mathfrak{NS}),$ we denote by $\mathrm{deg}p$ and $\overline{deg}(p)$ the degree of $p$ according as $\mathbb{Z}\cup Z$-gradation and  $\mathbb{Z}_2$-gradation respectively.
Let $V$ be a highest weight module  with highest weight vector $u$ of weight
$(c,h).$ Without loss of generality, we may assume that $u\in V_{\bar{0}}$ in the following.
We will
consider the tensor product modules $V\otimes SA_{a,b}$ and
$V\otimes SA_{a,b}^{'}$. Let us first introduce an auxiliary module, using the called ``shifting
technique" in [13].

\vskip 3mm\noindent{\bf{Lemma 3.1.}} The vector space $V\otimes
\mathbb{C}[t^{\pm\frac{1}{2}}]$ can be endowed with a $\mathfrak{NS}$-module
structure via
\begin{align*}
&L_k(pu\otimes t^s)=\left\{\begin{array}
{l@{\quad\quad}l} (L_k+a+s+kb-\mathrm{deg}p)pu\otimes t^{s+k}, & \mathrm{if}\ s+\mathrm{deg}p\in \mathbb{Z}, \\
(L_k+a+s+k(b-\frac{1}{2})-\mathrm{deg}p)pu\otimes t^{s+k}, &
\mathrm{if}\ s+\mathrm{deg}p\in Z,
\end{array}\right.\\
&G_{\frac{1}{2}+k}(pu\otimes t^s)=\left\{\begin{array}
{l@{\quad\quad}l} (G_{\frac{1}{2}+k}+(-1)^{\overline{deg}(p)})pu\otimes t^{s+k+\frac{1}{2}},\ \ \ \ \  \mathrm{if}\ s+\mathrm{deg}p\in \mathbb{Z}, \\
(G_{\frac{1}{2}+k}+(-1)^{\overline{deg}(p)}(a+s+(2k+1)(b-\frac{1}{2})-\mathrm{deg}p))pu\otimes t^{s+k+\frac{1}{2}},\\
\ \ \ \ \ \ \ \ \ \ \ \ \ \ \ \ \ \ \ \ \ \ \ \ \ \ \ \ \ \ \ \ \ \ \ \ \ \ \ \ \ \ \ \ \ \ \ \mathrm{if}\ s+\mathrm{deg}p\in
Z.
\end{array}\right.
\end{align*}

\vskip 3mm\noindent{\bf{Proof.}} It can be checked by straightforward but tedious calculations.
 \hfill$\Box$

\vskip 3mm\noindent{\bf{Lemma 3.2.}} The $\mathfrak{NS}$-module $V\otimes
\mathbb{C}[t^{\pm\frac{1}{2}}]$ is isomorphic to $V\otimes SA_{a,b}$ via the following map: for any $m\in \mathbb{Z},$
\begin{align*}
f:V\otimes SA_{a,b} &\rightarrow V\otimes \mathbb{C}[t^{\pm\frac{1}{2}}]\\
pu\otimes x_m & \mapsto pu\otimes t^{m+\mathrm{deg}p} \\
pu\otimes y_{\frac{1}{2}+m} & \mapsto pu\otimes t^{\frac{1}{2}+m+\mathrm{deg}p}.
\end{align*}

\vskip 3mm\noindent{\bf{Proof.}} It can be checked directly. \hfill$\Box$

When $a=b=\frac{1}{2},$ we see that $V\otimes \mathbb{C}[t^{\pm\frac{1}{2}}]$ has a submodule spanned by $pu\otimes t^{-\frac{1}{2}+\mathrm{deg}p}$
for all homogenous $p\in U(\mathfrak{n}^-),$ which is just a copy of $V.$ Let $\overline{V\otimes \mathbb{C}[t^{\pm\frac{1}{2}}]}$ be the corresponding quotient module in this case. For simplicity, we also denote by $\overline{V\otimes \mathbb{C}[t^{\pm\frac{1}{2}}]}$ when $V\otimes \mathbb{C}[t^{\pm\frac{1}{2}}]$ is irreducible.

The map in Lemma 3.2 induces an $\mathfrak{NS}$-module isomorphism of $V\otimes SA_{\frac{1}{2},\frac{1}{2}}^{'}$ and  $\overline{V\otimes \mathbb{C}[t^{\pm\frac{1}{2}}]}.$
Thus in the rest of this section we will consider $V\otimes \mathbb{C}[t^{\pm\frac{1}{2}}]$ and $\overline{V\otimes \mathbb{C}[t^{\pm\frac{1}{2}}]}$
instead of $V\otimes SA_{a,b}$ and $V\otimes SA_{a,b}^{'}.$ The advantage of the ``shifting technique" of the notation for $V\otimes \mathbb{C}[t^{\pm\frac{1}{2}}]$ is the weight space decomposition, that is
$$V\otimes t^s=\{x\in L|L_0x=(a+h+s)x\}.$$

\vskip 3mm From now on in this section, we fix $V=V(c,h),$  the simple highest weight module with highest weight vector $u$ of weight $(c,h).$

\vskip 3mm\noindent{\bf{Lemma 3.3.}} The $\mathfrak{NS}$-module $\overline{V(c,h)\otimes \mathbb{C}[t^{\pm\frac{1}{2}}]}$ is generated by $\{u\otimes t^{m}|m\in Z\}.$

\vskip 3mm\noindent{\bf{Proof.}} It can be proved by induction on $\mathrm{deg}p$. \hfill$\Box$

\vskip 3mm\noindent{\bf{Lemma 3.4.}} For any nonzero $\mathfrak{NS}$-submodule $W$ of $\overline{V(c,h)\otimes \mathbb{C}[t^{\pm\frac{1}{2}}]},$ there exists $N\in\mathbb{Z}\cup Z$ such that $V(c,h)\otimes t^s\subseteq W$ for all $s\in\mathbb{Z}\cup Z$ and $s\geq N.$

\vskip 3mm\noindent{\bf{Proof.}} By Lemma 3.2 and the well known result that the submodules of a weight module are also weight modules, there exist $W_s\subseteq V(c,h)$ for all $s\in \mathbb{Z}\cup Z$ such that $$W=\bigoplus_{s\in\mathbb{Z}\cup Z}W_s\otimes t^s.$$
For any nonzero vector $w\in W_s$ we can find homogenous $p_i\in U(\mathfrak{n}^{-})$ such that $$w=\sum_{i=1}^{r}p_iu=\sum_{i=1}^{r}w_i,$$
where $0\geq\mathrm{deg}p_1>\mathrm{deg}p_2>\cdots>\mathrm{deg}p_r$ and $w_i=p_iu.$
Choose a vector $w$ with $r$ being minimal among all nonzero elements in all $W_s,s\in \mathbb{Z}\cup Z.$ Now we prove that $r=1.$ Suppose this $w\in W_k.$
Denote $$l_i=-\mathrm{deg}p_i.$$
where, without loss of generality, we suppose that $\mathrm{deg}(w_i)+k\in \mathbb{Z}$ if $1\leq i\leq r_1$ and $\mathrm{deg}(w_i)+k\in Z$ if $r_1+1\leq i\leq r.$
If $a=b=\frac{1}{2},$ by the definition of $\overline{V(c,h)\otimes\mathbb{C}[t^{\pm\frac{1}{2}}]},$  we can further assume that $\mathrm{deg}p_i\neq \frac{1}{2}+k,$ i.e., $l_i\neq-\frac{1}{2}-k$ for all $i=1,\cdots,r.$
For $m,n\in \mathbb{Z}$ and $m,n>-\mathrm{deg}p_r,$ we have
\begin{align*}
&G_{\frac{1}{2}+m}G_{\frac{1}{2}+n}(w\otimes t^k)=\sum_{i=1}^{r_1}(a+k+\frac{1}{2}+n+(2m+1)(b-\frac{1}{2})+l_i)w_i\otimes t^{m+n+1+k}\\
&\ \ \ \ \ \ \ \ \ \ \ \ \ \ \ \ \ \ \ \ \ \ \ \ \ \ +\sum_{i=r_1+1}^{r}(a+k+(2n+1)(b-\frac{1}{2})+l_i)w_i\otimes t^{m+n+1+k},
\end{align*}and
\begin{align*}
&L_{m+n+1}(w\otimes t^k)=\sum_{i=1}^{r_1}(a+k+(m+n+1)b+l_i)w_i\otimes t^{m+n+1+k}\\
&\ \ \ \ \ \ \ \ \ \ \ \ \ \ \ \ \ \ \ \ \ +\sum_{i=r_1+1}^{r}(a+k+(m+n+1)(b-\frac{1}{2})+l_i)w_i\otimes t^{m+n+1+k}.
\end{align*}
If $r_1\geq1,$ for $1\leq j\leq r_1,$
we set
\begin{align*}
&A:=(a+k+(m+n+1)b+l_{j})G_{\frac{1}{2}+m}G_{\frac{1}{2}+n}(w\otimes t^k),\\
&B:=(a+k+\frac{1}{2}+n+(2m+1)(b-\frac{1}{2})+l_{j})L_{m+n+1}(w\otimes t^k).
\end{align*}Then by a calculation we get that
$$A-B=\sum_{i=1}^{r_1}(l_i-l_{j})(n-m)(b-1)w_i\otimes t^{m+n+1+k}+\sum_{i=r_1+1}^{r}f_iw_i\otimes t^{m+n+1+k},$$
where
\begin{align*}
&f_i=(n-m)\{l_j(b-\frac{1}{2})+l_i(b-1)+2(a+k)(b-\frac{3}{4})+2(m+n+1)(b-\frac{1}{2})\}.
\end{align*}
By the the minimality of $r,$  we have the following facts:
\begin{equation} \label{eq:1}
l_i-l_{j}=0\ \mathrm{for}\ 1\leq i\leq r_1\ \mathrm{since}\ b\neq1,
\end{equation}
\begin{equation} \label{eq:1}
f_i=0.
\end{equation}
By (3.1) we see that $r_1=1.$ By (3.2), if $r_1\neq r$ then $b=\frac{1}{2}$ and $a+k+l_i=0$ for $r_1+1\leq i\leq r.$ Since $k+l_i\in Z$
we see that $a=\frac{1}{2}$ by our convention about $a$ and $b$. So $l_i=-\frac{1}{2}-k,$ which is absurd. Thus  $r=1.$
Similarly, if $r-r_1\geq1,$ then we can also deduce that $r=1.$
Furthermore, $w=p_1u\in W_k$ such that $\mathrm{deg}p_1\neq \frac{1}{2}+k$ when $a=b=\frac{1}{2}.$

There exists $N_1\in \mathbb{Z}\cup Z$ such that $E_sp_1u=0$ and
$$ a+k+sb-\mathrm{deg}p_1, a+k+s(b-\frac{1}{2})-\mathrm{deg}p_1,  a+k+2s(b-\frac{1}{2})-\mathrm{deg}p_1$$ are all nonzero for any $s\geq N_1-k.$ Thus by
$$E_s(p_1u\otimes t^k)=\left\{\begin{array}
{l@{\quad\quad}l} (a+k+sb-\mathrm{deg}p_1)p_1u\otimes t^{s+k},\  \mathrm{if}\ s\in \mathbb{Z},\mathrm{deg}p_1\in \mathbb{Z}, \\
(a+k+s(b-\frac{1}{2})-\mathrm{deg}p_1)p_1u\otimes t^{s+k},\  \mathrm{if}\ s\in \mathbb{Z},\mathrm{deg}p_1\in Z, \\
(-1)^{\overline{deg}(p_1)}p_1u\otimes t^{s+k},\  \mathrm{if}\ s\in Z,\mathrm{deg}p_1\in Z, \\
(-1)^{\overline{deg}(p_1)}(a+k+(2s+1)(b-\frac{1}{2})-\mathrm{deg}p_1)p_1u\otimes t^{s+k},\  \mathrm{if}\ s\in Z,\mathrm{deg}p_1\in \mathbb{Z},
\end{array}\right.$$
we see that $E_s(p_1u\otimes t^k)$ is a nonzero element of $W$ and hence $p_1u\in W_s$ for all $s\geq N_1.$

For any $0<i\in \mathbb{Z}\cup Z$ and $s\geq N_1,$ we have $p_1u\in W_{s+i}$ and $E_i(p_1u\otimes t^s)\in W.$ So $E_ip_1u\in W_{s+i}$ for all $s\geq N_1.$
Inductively, we can show $xp_1u\in W_{s+\mathrm{deg}x}$ for all homogeneous $x\in U(\mathfrak{n}^{+})$ and $s\geq N_1.$ Since $V(c,h)$ is an irreducible $\mathfrak{NS}$-module, we can find some homogeneous $y\in U(\mathfrak{n}^{+})$ such that $yp_1u=u$
and hence $u\in W_s$ for all $s\geq N=N_1+\mathrm{deg}y.$

For any $0<i\in \mathbb{Z}\cup Z$ and $s\geq N,$ we have $E_{-i}(u\otimes t^{s+i})\in W,$ which means $E_{-i}u\in W_s$ for all $s\geq N.$ Proceeding by downward induction on $\mathrm{deg}p, p\in U(\mathfrak{n}^-),$ we can deduce that $pu\in W_s$ for $s\geq N$ and homogeneous $p\in U(\mathfrak{n}^-),$ thus $W_s=V(c,h)$ for all $s\geq N,$ as desired.
\hfill$\Box$

\vskip 3mm From the proof of Lemma 3.4 we can get the following corollary:

\vskip 3mm\noindent{\bf{Corollary 3.5.}} For any $0<s\in \mathbb{Z}\cup Z,$ the submodule generated by all $V(c,h)\otimes t^k, k>s$ is the same as the submodule generated by all $u\otimes t^k, k>s.$ We denote this submodule by $W^{(s)}.$

\vskip 3mm Now we introduce our linear map $\tilde{\varphi}_s: M(c,h)\rightarrow \mathbb{C}$ for $s\in\mathbb{Z}\cup Z.$ Denote by $T(n^-)$ the tensor algebra of $n^{-}.$ For any $s\in \mathbb{Z}\cup Z,$ define $\mathrm{deg}(E_s)=s$ and $\mathrm{deg}(1)=0,$ then $T(n^{-})$ becomes a $\mathbb{Z}\cup Z$-graded algebra. Now for any $s\in \mathbb{Z}\cup Z,\ a,b\in \mathbb{C},k\in \mathbb{Z}$ and homogeneous element $p\in T(n^{-}),$ we can inductively define a linear map $\varphi_s: T(n^{-})\rightarrow \mathbb{C}$ as following:
\begin{align*}
&\varphi_s(1)=1,\\
&\varphi_s(L_{-k}p)=\left\{\begin{array}
{l@{\quad\quad}l} -(a+s-kb-\mathrm{deg}p+k)\varphi_s(p),\  \mathrm{if}\ s+\mathrm{deg}p\in \mathbb{Z}, \\
-(a+s-k(b-\frac{1}{2})-\mathrm{deg}p+k)\varphi_s(p),\  \mathrm{if}\ s+\mathrm{deg}p\in Z,
\end{array}\right.\\
&\varphi_s(G_{-\frac{1}{2}-k}p)=\left\{\begin{array}
{l@{\quad\quad}l} -(-1)^{\overline{deg}(p)}\varphi_s(p),\  \mathrm{if}\ s+\mathrm{deg}p\in \mathbb{Z}, \\
-(-1)^{\overline{deg}(p)}(a+s-(2k+1)(b-\frac{1}{2})-\mathrm{deg}p+k)\varphi_s(p),\  \mathrm{if}\ s+\mathrm{deg}p\in Z,
\end{array}\right.
\end{align*}
Denote by $J$ the two-sided ideal generated by the following elements:
\begin{align*}
&L_{-m}L_{-n}-L_{-n}L_{-m}-(m-n)L_{-m-n},\\
&L_{-m}G_{-r}-G_{-r}L_{-m}-(\frac{m}{2}-r)G_{-m-r},\\
&G_{-r}G_{-s}+G_{-s}G_{-r}-2L_{-r-s},
\end{align*} where $m,n\in \mathbb{N}$ and $r,s\in Z_{>0},$
then $U(\mathfrak{n}^-)$ is just the quotient $T(\mathfrak{n}^-)/J.$
By the definition of $\varphi_s(s\in \mathbb{Z}\cup Z),$ for any homogeneous element $p\in T(\mathfrak{n}^-),$ we can verify the following identities:
\begin{align*}
&\varphi_s((L_{-m}L_{-n}-L_{-n}L_{-m})p)=(m-n)\varphi_s(L_{-m-n}p),\\
&\varphi_s(L_{-m}G_{-r}-G_{-r}L_{-m})p=(\frac{m}{2}-r)\varphi_s(G_{-m-r}p),\\
&\varphi_s(G_{-r}G_{-s}+G_{-s}G_{-r}p)=2\varphi_s(L_{-r-s}p).
\end{align*}
Thus $J\subseteq\mathrm{ker}\varphi_s$ and $\varphi_s$ induces a linear map on $U(\mathfrak{n}^-)$, denoted also by $\varphi_s.$
Since the Verma module $M(c,h)$ is free of rank $1$ as a $U(\mathfrak{n}^{-})$-module, we can define the linear map
$$\tilde{\varphi}_s: M(c,h)\rightarrow \mathbb{C}, \tilde{\varphi}_s(pu)=\varphi_s(p), \forall p\in U(\mathfrak{n}^{-}).$$

\vskip 3mm\noindent{\bf{Lemma 3.6.}} Let $W$ be a $NS$-submodule of $\overline{V(c,h)\otimes \mathbb{C}[t^{\pm\frac{1}{2}}]}\cong V(c,h)\otimes SA_{a,b}^{'}$ with
$W\supseteq W^{(s)}$ for some $0<s\in\mathbb{Z}\cup Z.$ Then $$pu\otimes t^s\equiv \varphi_s(p)u\otimes t^{s}(\mathrm{mod}\ W), \forall p\in U(\mathfrak{n}^-).$$

\vskip 3mm\noindent{\bf{Proof.}} By the definition of $\varphi_s,$ it is obvious that $u\otimes t^s\equiv \varphi_s(1)u\otimes t^{s}(\mathrm{mod}\ W)$.
Let $q\in U(\mathfrak{n}^-)$ be a homogeneous element and suppose that the result holds for this $q$. Then for any $0<k\in\mathbb{Z},$ we have
\begin{align*}
&L_{-k}(qu\otimes t^{s+k})=\left\{\begin{array}
{l@{\quad\quad}l} (L_{-k}+a+s+k-kb-\mathrm{deg}q)qu\otimes t^{s}, & \mathrm{if}\ s+\mathrm{deg}q\in \mathbb{Z}, \\
(L_{-k}+a+s+k-k(b-\frac{1}{2})-\mathrm{deg}q)qu\otimes t^{s}, &
\mathrm{if}\ s+\mathrm{deg}q\in Z,
\end{array}\right.\end{align*}
which means $L_{-k}(qu\otimes t^{s+k})\equiv \varphi_s(L_{-k}q)u\otimes t^{s+k}(\mathrm{mod}\ W).$
For any $0<k\in Z,$ we have
\begin{align*}
&G_{-k}(qu\otimes t^{s+k})=\left\{\begin{array}
{l@{\quad\quad}l} (G_{-k}+(-1)^{\overline{deg}(q)})qu\otimes t^{s},\ \ \ \ \ \ \ \ \ \ \ \ \ \ \ \ \  \ \ \ \ \ \mathrm{if}\ s+\mathrm{deg}q\in \mathbb{Z}, \\
(G_{-k}+(-1)^{\overline{deg}(q)}(a+s+k-2k(b-\frac{1}{2})-\mathrm{deg}q))qu\otimes t^{s},\\
\ \ \ \ \ \ \ \ \ \ \ \ \ \ \ \ \ \ \ \ \ \ \ \ \ \ \ \ \ \ \ \ \  \ \ \ \ \ \ \ \ \ \ \ \ \ \ \ \ \ \   \  \ \ \ \ \mathrm{if}\ s+\mathrm{deg}q\in
Z.
\end{array}\right.
\end{align*}
which means $G_{-k}(qu\otimes t^{s+k})\equiv \varphi_s(G_{-k}q)u\otimes t^{s+k}(\mathrm{mod}\ W).$
So the lemma follows from induction on $\mathrm{deg}q$.
\hfill$\Box$

Let $J(c,h)$ be the maximal submodule of $M(c,h).$ By [3], $J(c,h)$ is generated by at most two homogeneous singular vectors $Q_1, Q_2\in U(\mathfrak{n}^-).$ If  $J(c,h)$ is generated by exactly one homogeneous singular vectors, we assume that $Q_1=Q_2.$ If $M(c,h)$ is simple, then we assume that $Q_1=Q_2=0.$ By our assumptions, $Q_1$ and $Q_2$ are unique up to nonzero scalars.
By the definition of $\varphi_s$,  there exists only finitely many $s\in \mathbb{Z}\cup Z$ such that $\varphi_s(Q_1)=\varphi_s(Q_2)=0.$

Now we give the main result of this section.

\vskip 3mm\noindent{\bf{Theorem 3.7.}} Let $a,b,c,h\in \mathbb{C}$ with $0\leq \mathrm{Re}(a)<1$ and $b\neq 1.$ Then we have

(a) $V(c,h)\otimes SA_{a,b}^{'}$ is simple if and only if there is no (other than $-\frac{1}{2}$ if $(a,b)=(\frac{1}{2},\frac{1}{2})$) number $s\in\mathbb{Z}\cup Z$ such that $\varphi_s(Q_1)=\varphi_s(Q_2)=0.$

(b) If $s$ is the maximal number (not equal to $-\frac{1}{2}$ if $(a,b)=(\frac{1}{2},\frac{1}{2})$) in $\mathbb{Z}\cup Z$ such that $\varphi_s(Q_1)=\varphi_s(Q_2)=0,$ then $W^{(s)}$ is the unique simple submodule and all its nonzero weight space are infinite-dimensional.

\vskip 3mm\noindent{\bf{Proof.}} Let $W$ be a nonzero $NS$-submodule of $\overline{V(c,h)\otimes \mathbb{C}[t^{\pm\frac{1}{2}}]}\cong V(c,h)\otimes SA_{a,b}^{'}.$ By Lemma 3.4, there exists $N\in \mathbb{Z}\cup Z$ such that $W^{(N)}\subseteq W.$  Then by Lemma 3.6 we see that
$$Q_iu\otimes t^N-\varphi_N(Q_i)u\otimes t^N=-\varphi_N(Q_i)u\otimes t^N\in W(i=1,2)$$
since $V(c,h)=M(c,h)/J(c,h)$ and $Q_iu\in J(c,h).$
If $(a,b)\neq(\frac{1}{2},\frac{1}{2})$ and there is no $s\in\mathbb{Z}\cup Z$ such that $\varphi_s(Q_1)=\varphi_s(Q_2)=0,$ then we have $u\otimes t^N\in W,$  i.e., $W^{(N-1)}\subseteq W.$ By induction we deduce that
$u\otimes t^s\in W$ for all $s\in\mathbb{Z}\cup Z.$ If $(a,b)=(\frac{1}{2},\frac{1}{2})$ and and there is no $-\frac{1}{2}\neq s\in\mathbb{Z}\cup Z$ such that $\varphi_s(Q_1)=\varphi_s(Q_2)=0,$ we can also get that $u\otimes t^s\in W$ for all $s\in\mathbb{Z}\cup Z$ since $u\otimes t^{-\frac{1}{2}}=0\in W$ by the definition of $\overline{V(c,h)\otimes\mathbb{C}[t^{\pm\frac{1}{2}}]}$.
By Corollary 3.5, we see that $W=\overline{V(c,h)\otimes \mathbb{C}[t^{\pm\frac{1}{2}}]}$ and $V(c,h)\otimes SA_{a,b}^{'}$ is irreducible.

Let $N\in \mathbb{Z}\cup Z$ be any number such that $\varphi_N(Q_1)=\varphi_N(Q_2)=0.$ We will show that $W^{(N)}$ is a proper submodule of $\overline{V(c,h)\otimes\mathbb{C}[t^{\pm\frac{1}{2}}]}.$ By the definitions of $\varphi_N$ and $\tilde{\varphi}_N$, we have
$$\tilde{\varphi}_N(J(c,h))=\tilde{\varphi}_N(U(\mathfrak{n}^{-})Q_1u+U(\mathfrak{n}^{-})Q_2u)=\varphi_N(U(\mathfrak{n}^{-})Q_1)
+\varphi_N(U(\mathfrak{n}^{-})Q_2)=0.$$
Thus $\tilde{\varphi}_N$ induces a linear map $\tilde{\varphi}_N: V(c,h)\rightarrow\mathbb{C}$ which sends $pu$ to $\varphi_N(p)$ for any $p\in U(\mathfrak{n}^-),$ moreover,
$\mathrm{ker}(\tilde{\varphi}_N)\neq V(c,h)$ since $\tilde{\varphi}_N(u)=1.$

By the PBW theorem, the weight space of $W^{(N)}$ with weight $a+h+N$ is
\begin{align*}
W^{(N)}_N\otimes t^{N}&=\sum_{0<s\in\mathbb{Z}\cup Z}E_{-s}(V(c,h)\otimes t^{N+s})\\
&=\sum_{0<s\in \mathbb{Z}}\sum_{\begin{subarray}\
p\in U(\mathfrak{n}^{-})\\
N+\mathrm{deg}p\in \mathbb{Z}
\end{subarray}}((L_{-s}+a+N+s-\mathrm{deg}p-sb)pu\otimes t^{N})\\
&\ \ +\sum_{0<s\in \mathbb{Z}}\sum_{\begin{subarray}\
p\in U(\mathfrak{n}^{-})\\
N+\mathrm{deg}p\in Z
\end{subarray}}((L_{-s}+a+N+s-\mathrm{deg}p-s(b-\frac{1}{2}))pu\otimes t^{N})\\
&\ \ +\sum_{0<s\in Z}\sum_{\begin{subarray}\
p\in U(\mathfrak{n}^{-})\\
N+\mathrm{deg}p\in Z
\end{subarray}}((G_{-s}+1)pu\otimes t^{N})\\
&\ \ +\sum_{0<s\in Z}\sum_{\begin{subarray}\
p\in U(\mathfrak{n}^{-})\\
N+\mathrm{deg}p\in \mathbb{Z}
\end{subarray}}((G_{-s}-(a+N+s-2s(b-\frac{1}{2})-\mathrm{deg}p))pu\otimes t^{N}).
\end{align*}
It is easy to see that if $N+\mathrm{deg}p\in \mathbb{Z},$ then
$$\tilde{\varphi}_N(L_{-s}+a+N+s-\mathrm{deg}p-sb)pu)=\varphi_N(L_{-s}p)+(a+N+s-\mathrm{deg}p-sb)\varphi_N(p)=0,$$ and
$$\tilde{\varphi}_N((G_{-s}-(a+N+s-2s(b-\frac{1}{2})-\mathrm{deg}pu)=0.$$
 Similarly,
if $N+\mathrm{deg}p\in Z,$ then
$$\tilde{\varphi}_N(L_{-s}+a+N+s-\mathrm{deg}p-s(b-\frac{1}{2}))pu)=0$$ and
$$\tilde{\varphi}_N(G_{-s}+1)pu)=0.$$
Thus $W^{(N)}_N\subseteq \mathrm{ker}(\tilde{\varphi}_N)\subsetneqq V(c,h),$ i.e., $W^{(N)}$ is a proper submodule of $\overline{V(c,h)\otimes\mathbb{C}[t^{\pm\frac{1}{2}}]}$.

The proof of (b) is similar as the proof of sufficiency of (a). In Fact, let $W$ be an arbitrary Proper submodule of $\overline{V(c,h)\otimes\mathbb{C}[t^{\pm\frac{1}{2}}]}.$ If $(a,b)\neq(\frac{1}{2},\frac{1}{2}),$ Let $N_0$ be the largest number such that
$\varphi_{N_0}(Q_1)=\varphi_{N_0}(Q_2)=0.$ By Lemma 3.4 there exists $N\in\mathbb{Z}\cup Z$ such that $V(c,h)\otimes t^s\subseteq W$ for all $s\in\mathbb{Z}\cup Z$ and $s>N.$ Then by Lemma 3.6 we have $$Q_iu\otimes t^{N}-\varphi_N(Q_i)u\otimes t^N\in W,\ i=1,2,$$ i.e., $u\otimes t^N\in W.$ By Corollary 3.5 $W\supseteq W^{(N-1)}.$ By induction, we see that $W\supseteq W^{(N_0)}.$ If $(a,b)=(\frac{1}{2},\frac{1}{2}),$ Let $N_0\neq -\frac{1}{2}$ be the largest number of $\mathbb{Z}\cup Z$ such that $\varphi_{N_0}(Q_1)=\varphi_{N_0}(Q_2)=0.$ Since $u\otimes t^{(-\frac{1}{2})}=0\in W,$ we also get by induction that $W\supseteq W^{(N_0)}.$ Thus $W^{(N_0)}$ is the smallest submodule of $\overline{V(c,h)\otimes\mathbb{C}[t^{\pm\frac{1}{2}}]}$ which has to be simple. From [8], we know that any weight spaces are infinite-dimensional.  \hfill$\Box$

\vskip 3mm Let $\Phi$ be the set of $s\in\mathbb{Z}\cup Z$ ($\neq-\frac{1}{2}$ if $(a,b)=(\frac{1}{2},\frac{1}{2})$) such that $\varphi_s(Q_1)=\varphi_s(Q_1)=0.$ We have the following remark:

\vskip 3mm\noindent{\bf{Remark 3.8.}}  We have a consequence of submodules of $V\otimes SA_{a,b}^{'}:$
$$\cdots\subseteq W^{(s+\frac{1}{2})}\subseteq W^{(s)}\subseteq W^{(s-\frac{1}{2})}\subseteq\cdots.$$

If $Q_1=Q_2=0,$ i.e., $V(c,h)=M(c,h),$  then $\Phi=\mathbb{Z}\cup Z$ ($\Phi=\mathbb{Z}\setminus \{-\frac{1}{2}\}$ if $(a,b)=(\frac{1}{2},\frac{1}{2})$)
and all inclusions $W^{(s+1)}\subseteq W^{(s)}$ are proper for $s\in\Phi.$
It is easy to see that $W^{(s)}/W^{(s+\frac{1}{2})}$ is a highest weight module of highest weight $a+h+s$ for any $s\in\Phi.$ In particular, $V\otimes SA_{a,b}^{'}$ is an extension of countably many highest weight modules.

If $Q_1\neq 0$ or $Q_2\neq0,$ then $\Phi$ is a finite subset of $\mathbb{Z}\cup Z$ and the only proper inclusions in the
above consequence of submodules are $W^{(s)}\subsetneqq W^{(s-\frac{1}{2})}$ for $s\in \Phi.$  If we write $\Phi=\{s_1,s_2,\cdots,s_r\}$ with $s_1<s_2<\cdots<s_r$ and take any $s_0<s_1,$ then the consequence of submodules can be simplified as:
$$0\subsetneqq W^{(s_r)}\subsetneqq\cdots\subsetneqq W^{(s_{2})}\subsetneqq W^{(s_{1})}\subsetneqq W^{(s_{0})}=V\otimes SA_{a,b}^{'}.$$

By Theorem 3.7 (b) we know that $W^{(n_r)}$ is the unique minimal nonzero proper submodule of $V\otimes SA_{a,b}^{'}.$

\vskip 3mm\noindent{\section{Isomorphism objects}}

\vskip 3mm In this section, we will determine the necessary and sufficient conditions for two $\mathfrak{NS}$-tensor modules with form $V(c,h)\otimes SA_{a,b}^{'}$ to be isomorphism.  The main result is the following theorem:

\vskip 3mm\noindent{\bf{Theorem 4.1.}} Let $c, h, a, b, c_1, h_1, a_1, b_1\in\mathbb{C}$ with $0\leq \mathfrak{Re}(a),\mathfrak{Re}(a_1)<1, b, b_1\neq1.$  Then $V(c,h)\otimes SA_{a,b}^{'}\cong V_1(c_1,h_1)\otimes SA_{a_1,b_1}^{'}$ if and only if $c=c_1,h=h_1,a=a_1$ and $b=b_1.$

As a direct consequence we have

\vskip 3mm\noindent{\bf{Theorem 4.2.}} Let $V(c,h)$ and $V_1(c_1,h_1)$ be any two highest weight modules (not necessarily simple) of the Neveu-Schwarz algebra,
and let $a,b,a_1,b_1\in\mathbb{C}.$ Then $V(c,h)\otimes SA_{a,b}^{'}\cong V_1(c_1,h_1)\otimes SA_{a_1,b_1}^{'}$ if and only if $V(c,h)\cong V_1(c_1,h_1)$ and $SA_{a,b}^{'}\cong SA_{a_1,b_1}^{'}.$

\vskip 3mm\noindent{\bf{Proof\ of\ Theorem 4.1.}} The sufficiency is trivial.
Assume $$\psi:V(c,h)\otimes SA_{a,b}^{'}\rightarrow V_1(c_1,h_1)\otimes SA_{a_1,b_1}^{'}$$ is an isomorphism of $\mathfrak{NS}$-modules.
 It is clear that
\begin{equation} \label{eq:1}
c=c_1.
\end{equation}
 Fix any $k\in\mathbb{Z}\cup Z$ such that $k\neq -\frac{1}{2}$ when $(a,b)=(\frac{1}{2},\frac{1}{2}).$ Since $\psi(u\otimes t^{k})$ and $u\otimes t^{k}$ are of the same weight, we can assume that
\begin{equation} \label{eq:1}
\psi(u\otimes t^{k})=\sum_{i=1}^{r}p_{i,k}u_1\otimes t^l,
\end{equation}
where $p_{i,k}$ are homogeneous elements of $U(\mathfrak{n}^{-})$ with $0\geq\mathrm{deg}p_{1,k}>\mathrm{deg}p_{2,k}>\cdots>\mathrm{deg}p_{r,k},$ and $l$ satisfies
\begin{equation} \label{eq:1}
a+h+k=a_1+h_1+l.
\end{equation}
Without loss of generality, we suppose that $\mathrm{deg}p_{i,k}+l\in \mathbb{Z}$ if $1\leq i\leq r_1$ and $\mathrm{deg}p_{i,k}+l\in Z$ if $r_1+1\leq i\leq r.$
By the definition of $\overline{V(c,h)\otimes\mathbb{C}[t^{\pm\frac{1}{2}}]},$  we can further assume that
\begin{equation} \label{eq:1}
\mathrm{deg}p_{i,k}\neq \frac{1}{2}+l\ \ \mathrm{if}\ a=b=\frac{1}{2}.
\end{equation}
For $m,n\in \mathbb{Z}$ and $m,n>-\mathrm{deg}p_r,$ we have
\begin{align*}
&\psi(G_{\frac{1}{2}+m}G_{\frac{1}{2}+n}(u\otimes t^k))\\
=&\left\{\begin{array}
{l@{\quad\quad}l} (a+k+\frac{1}{2}+n+(2m+1)(b-\frac{1}{2}))\psi(u\otimes t^{m+n+1+k}), &\mathrm{if}\ k\in \mathbb{Z}, \\
(a+k+(2n+1)(b-\frac{1}{2}))\psi(u\otimes t^{m+n+1+k}), & \mathrm{if}\ k\in Z
\end{array}\right.\\
=&G_{\frac{1}{2}+m}G_{\frac{1}{2}+n}\psi(u\otimes t^k))\\
=&\sum_{i=1}^{r_1}(a_1+l+\frac{1}{2}+n+(2m+1)(b_1-\frac{1}{2})-\mathrm{deg}p_{i,k})w_i\otimes t^{m+n+1+l}
\end{align*}
\begin{equation} \label{eq:1}
+\sum_{i=r_1+1}^{r}(a_1+l+(2n+1)(b_1-\frac{1}{2})-\mathrm{deg}p_{i,k})w_i\otimes t^{m+n+1+l},
\end{equation}
and
\begin{align*}
&\psi(L_{m+n+1}(u\otimes t^k))\\
=&\left\{\begin{array}
{l@{\quad\quad}l} (a+k+(m+n+1)b)\psi(u\otimes t^{m+n+1+k}), & \mathrm{if}\ k\in \mathbb{Z}, \\
(a+k+(m+n+1)(b-\frac{1}{2}))\psi (u\otimes t^{m+n+1+k}), &
\mathrm{if}\ k\in Z
\end{array}\right.\\
=&L_{m+n+1}\psi(u\otimes t^k))\\
=&\sum_{i=1}^{r_1}(a_1+l+(m+n+1)b_1-\mathrm{deg}p_{i,k})w_i\otimes t^{m+n+1+k}
\end{align*}
\begin{equation} \label{eq:1}
 +\sum_{i=r_1+1}^{r}(a_1+l+(m+n+1)(b_1-\frac{1}{2})-\mathrm{deg}p_{i,k})w_i\otimes t^{m+n+1+k}.
\end{equation}

\vskip 3mm\noindent{\bf{Claim 1.}} $r=r_1$ when $k\in\mathbb{Z}$ ;\ \ $r_1=0$ when $k\in Z.$
We prove the first part. The proof for the second part is similar.

Otherwise, $r>r_1.$ For any $r_1+1\leq i\leq r,$  by (4.5) and (4.6), we have
\begin{align*}
&(a+k+(m+n+1)b)(a_1+l+(2n+1)(b_1-\frac{1}{2})-\mathrm{deg}p_{i,k}))\\
=&(a+k+\frac{1}{2}+n+(2m+1)(b-\frac{1}{2}))(a_1+l+(m+n+1)(b_1-\frac{1}{2})-\mathrm{deg}p_{i,k}).
\end{align*}
We can rewrite it as following:
$$2(b-\frac{1}{2})(b_1-\frac{1}{2})(n^2-m^2)+((b-1)(a_1+l-\mathrm{deg}p_{i,k})+(a+k+2b-1)(b_1-\frac{1}{2}))(n-m)=0.$$
Thus we have
$(b-\frac{1}{2})(b_1-\frac{1}{2})=0$ and
$(b-1)(a_1+l-\mathrm{deg}p_{i,k})+(a+k+2b-1)(b_1-\frac{1}{2})=0.$

If $b_1=\frac{1}{2},$ then we have $a_1=\frac{1}{2}$ since $b\neq1,\ l-\mathrm{deg}p_{i,k}\in Z$ and $0\leq\mathfrak{Re}a_1<1.$
Thus $l=-\frac{1}{2}+\mathrm{deg}p_{i,k}$ which is contrary to (4.3).

If $b_1\neq\frac{1}{2},$ then $b=\frac{1}{2}$ and
\begin{equation} \label{eq:1}
(a+k)(b_1-\frac{1}{2})=\frac{1}{2}(a_1+l-\mathrm{deg}p_{i,k}).
\end{equation}
In order to derive a contradiction, we consider $\psi(L_mL_n(u\otimes t^k))$ and $\psi(L_{m+n}(u\otimes t^k)).$ We have
\begin{align*}
&\psi(L_{m}L_n(u\otimes t^k))\\
=&(a+n+k+mb)(a+k+nb)\psi(u\otimes t^{m+n+k})\\
=&L_{m}L_n(\psi(u\otimes t^k))\\
=&\sum_{i=1}^{r_1}(a_1+n+l+mb_1-\mathrm{deg}p_{i,k})(a_1+l+nb_1-\mathrm{deg}p_{i,k})w_i\otimes t^{m+n+k}\\
&+\sum_{i=r_1+1}^{r}(a_1+n+l+m(b_1-\frac{1}{2})-\mathrm{deg}p_{i,k})(a_1+l+n(b_1-\frac{1}{2})-\mathrm{deg}p_{i,k})w_i\otimes t^{m+n+k},
\end{align*}
and
\begin{align*}
&\psi(L_{m+n}(u\otimes t^k))\\
=&(a+k+(m+n)b)\psi(u\otimes t^{m+n+k})\\
=&L_{m+n}(\psi(u\otimes t^k))\\
=&\sum_{i=1}^{r_1}(a_1+l+(m+n)b_1-\mathrm{deg}p_{i,k})w_i\otimes t^{m+n+k}\\
&+\sum_{i=r_1+1}^{r}(a_1+l+(m+n)(b_1-\frac{1}{2})-\mathrm{deg}p_{i,k})w_i\otimes t^{m+n+k},
\end{align*}where $b=\frac{1}{2}.$
So we have the following identity:
\begin{align*}
(a+n+k+\frac{1}{2}m)(a+k+\frac{1}{2}n)(a_1+l+(m+n)(b_1-\frac{1}{2})-\mathrm{deg}p_{i,k})
\end{align*}
\begin{equation} \label{eq:1}
 =(a+k+\frac{1}{2}(m+n))(a_1+n+l+m(b_1-\frac{1}{2})-\mathrm{deg}p_{i,k})(a_1+l+n(b_1-\frac{1}{2})-\mathrm{deg}p_{i,k}).
\end{equation}
If $a+k=0,$ by (4.7), then $a_1+l-\mathrm{deg}p_{i,k}=0$ and (4.8) becomes
$$n+\frac{1}{2}m=n+m(b_1-\frac{1}{2}).$$
Which implies that $b_1=1,$ a contradiction. Thus $a+k\neq0$ and $a_1+l-\mathrm{deg}p_{i,k}\neq0.$ We can view (4.8) as a equation about variables $m$ and $n,$ by comparing the constant term on two sides we get that $$(a+k)^{2}(a_1+l-\mathrm{deg}p_{i,k})=(a+k)(a_1+l-\mathrm{deg}p_{i,k})^{2}.$$
which implies that
\begin{equation} \label{eq:1}
a+k=a_1+l-\mathrm{deg}p_{i,k}\neq0.
\end{equation}
Then we get that $b_1=1$ by (4.7) and (4.9), which is also a contradiction. Thus $r=r_1.$ The proof of the claim for the case $k\in Z$ is similar, we omit the details.

By Claim 1, (4.5) and (4.6) are rewritten as following
\begin{eqnarray}
&&\psi(G_{\frac{1}{2}+m}G_{\frac{1}{2}+n}(u\otimes t^k))\nonumber\\
&&=\left\{\begin{array}
{l@{\quad\quad}l} (a+k+\frac{1}{2}+n+(2m+1)(b-\frac{1}{2}))\psi(u\otimes t^{m+n+1+k}), &\mathrm{if}\ k\in \mathbb{Z}, \\
(a+k+(2n+1)(b-\frac{1}{2}))\psi(u\otimes t^{m+n+1+k}), & \mathrm{if}\ k\in Z
\end{array}\right.\nonumber\\
&&=G_{\frac{1}{2}+m}G_{\frac{1}{2}+n}(\psi(u\otimes t^k))\nonumber\\
&&=\left\{\begin{array}
{l@{\quad\quad}l}\sum_{i}(a_1+l+\frac{1}{2}+n+(2m+1)(b_1-\frac{1}{2})-\mathrm{deg}p_{i,k})p_{i,k}u_1\otimes t^{m+n+1+l}, \mathrm{if}\ k\in \mathbb{Z}, \\
\sum_{i}(a_1+l+(2n+1)(b_1-\frac{1}{2})-\mathrm{deg}p_{i,k})p_{i,k}u_1\otimes t^{m+n+1+l}, \mathrm{if}\ k\in Z
\end{array}\right.\label{1}
\end{eqnarray}
and
\begin{eqnarray}
&&\psi(L_{m+n+1}(u\otimes t^k))\nonumber\\
&&=\left\{\begin{array}
{l@{\quad\quad}l} (a+k+(m+n+1)b)\psi(u\otimes t^{m+n+1+k}), & \mathrm{if}\ k\in \mathbb{Z}, \\
(a+k+(m+n+1)(b-\frac{1}{2}))\psi (u\otimes t^{m+n+1+k}), &
\mathrm{if}\ k\in Z
\end{array}\right.\nonumber\\
&&=L_{m+n+1}(\psi(u\otimes t^k))\nonumber\\
&&=\left\{\begin{array}
{l@{\quad\quad}l} \sum_{i}(a_1+l+(m+n+1)b_1-\mathrm{deg}p_{i,k})p_{i,k}u_1\otimes t^{m+n+1+k}, \mathrm{if}\ k\in \mathbb{Z}, \\
\sum_{i}(a_1+l+(m+n+1)(b_1-\frac{1}{2})-\mathrm{deg}p_{i,k})p_{i,k}u_1\otimes t^{m+n+1+k},  \mathrm{if}\ k\in Z
\end{array}\right.\label{1}
\end{eqnarray}
where $l+\mathrm{deg}p_{i,k}\in\mathbb{Z}$ if $k\in\mathbb{Z}$ and $l+\mathrm{deg}p_{i,k}\in Z$ if $k\in Z.$

\vskip 3mm\noindent{\bf{Claim 2.}} $b=b_1$ and $\psi(u\otimes t^k)=p_{1,k}u_1\otimes t^l.$ By (4.10) and (4.11),
for $k\in\mathbb{Z}$ we have
\begin{align*}
&(a+k+(m+n+1)b)(a_1+l+\frac{1}{2}+n+(2m+1)(b_1-\frac{1}{2})-\mathrm{deg}p_{i,k})\\
=&(a+k+\frac{1}{2}+n+(2m+1)(b-\frac{1}{2}))(a_1+l+(m+n+1)b_1-\mathrm{deg}p_{i,k}).
\end{align*}
We rewrite it as following:
$$(b-b_1)(n^2-m^2)+((a+k)(1-b_1)+(b-1)(a_1+l-\mathrm{deg}p_{i,k})+b-b_1)(m+n)=0.$$
Hence we get that
\begin{equation} \label{eq:1}
b=b_1,
\end{equation}
\begin{equation} \label{eq:1}
(a+k)(1-b_1)+(b-1)(a_1+l-\mathrm{deg}p_{i,k})=0,
\end{equation}
By (4.3), (4.12), (4.13) and the fact that $b,b_1\neq1,$ we get that
\begin{equation} \label{eq:1}
h_1-h+\mathrm{deg}p_{i,k}=0,
\end{equation}
which implies that $r=1$, i.e.,
$\psi(u\otimes t^k)=p_{1,k}u_1\otimes t^l$ for $k\in\mathbb{Z}.$
It can be proved similarly that  $\psi(u\otimes t^k)=p_{1,k}u_1\otimes t^l$ for $k\in Z.$ Thus Claim 2 holds.

\vskip 3mm\noindent{\bf{Claim 3.}} $p_{1,k}$ is a nonzero scalar for any $k\in\mathbb{Z}\cup Z$.
Note that $\overline{V(c,h)\otimes\mathbb{C}[t^{\pm\frac{1}{2}}]}\cong V\otimes SA_{a,b}^{'}$ is generated by $\{u\otimes t^k|k\in\mathbb{Z}\cup Z\}$,
we have
$$\sum_{k\in\mathbb{Z}\cup Z}U(\mathfrak{n}^{-})(u\otimes t^{k})=V(c,h)\otimes SA_{a,b}^{'}\ \mathrm{for}\ (a,b)\neq(\frac{1}{2},\frac{1}{2})$$
and
$$\sum_{k\in\mathbb{Z}\cup Z\setminus\{-\frac{1}{2}\}}U(\mathfrak{n}^{-})(u\otimes t^{k})=V(c,h)\otimes SA_{a,b}^{'}\ \mathrm{for}\ (a,b)=(\frac{1}{2},\frac{1}{2}),$$
which together with Claim 2, we have
$$\sum_{k\in\mathbb{Z}\cup Z}U(\mathfrak{n}^{-})p_{1,k}u_1\otimes t^{l}=V(c,h)\otimes SA_{a,b}^{'}\ \mathrm{for}\ (a,b)\neq(\frac{1}{2},\frac{1}{2})$$
and
$$\sum_{k\in\mathbb{Z}\cup Z\setminus\{-\frac{1}{2}\}}U(\mathfrak{n}^{-})p_{1,k}u_1\otimes t^{l}=V(c,h)\otimes SA_{a,b}^{'}\ \mathrm{for}\ (a,b)=(\frac{1}{2},\frac{1}{2}).$$
Thus $p_{1,k}$ is a nonzero scalar for any $k\in\mathbb{Z}\cup Z\ (k\neq-\frac{1}{2}$ if $(a,b)=(\frac{1}{2},\frac{1}{2})).$
Thus Claim 3 holds.

By (4.12), (4.13), (4.14) and Claim 3, we have
\begin{equation} \label{eq:1}
a+k=a_1+l,
\end{equation}and
\begin{equation} \label{eq:1}
h_1=h.
\end{equation}
Recall that $l+\mathrm{deg}p_{1,k}=l\in \mathbb{Z}$(resp., $Z$) if $k\in\mathbb{Z}$(resp., $Z$), $0\leq\mathfrak{Re}a<1,$ $0\leq\mathfrak{Re}a_1<1,$  Claim 3 and (4.15) imply that
\begin{equation} \label{eq:1}
a=a_1.
\end{equation}
Thus the necessity of Theorem 4.1 follows from (4.1), (4.12), (4.16) and (4.17).
\hfill$\Box$

\end{document}